\documentclass[10pt]{amsart}
\usepackage{graphicx}
\usepackage{amsmath}
\usepackage{amsthm}
\usepackage{amsfonts}
\usepackage{amssymb}
\usepackage{cite}

\usepackage{epic}
\usepackage{pstricks}

\theoremstyle{plain}
\newtheorem{theorem}{Theorem}[section]

\newtheorem{lemma}[theorem]{Lemma}

\newtheorem{definition}[theorem]{Definition}
\newtheorem{remark}[theorem]{Remark}

\newcommand{\lp}{\left(}
\newcommand{\rp}{\right)}

\newcommand{\nn}{\overrightarrow{n}}

\renewcommand{\v}[2]{\overrightarrow{V_{#1}V_{#2}}}

\newcommand{\intr}[2]{\overline{#1,#2}}

\newcommand{\s}{\,\rule[-3.5pt]{1pt}{12pt}\,}

\newcommand{\De}{\Delta}
\newcommand{\al}{\alpha}
\newcommand{\be}{\beta}

\newcommand{\si}{\sigma}

\newcommand{\g}{\gamma}
\newcommand{\vp}{\varepsilon}
\renewcommand{\th}{\theta}

\renewcommand{\le}{\leqslant}
\renewcommand{\ge}{\geqslant}

\newcommand{\R}{\mathbb{R}}
\newcommand{\Z}{\mathbb{Z}}

\renewcommand{\P}{\mathcal{P}}

\newcommand{\conv}{\operatorname{conv}}
\renewcommand{\int}{\operatorname{int}}

\newcommand{\sign}{\operatorname{sign}}

\begin{document}


\title[Convex Cyclic Polygons: Angles]
{A Characterization of the Convexity of Cyclic Polygons in Terms of the 
Central Angles} 


\author{Iosif Pinelis}  
                   

\date{\today; file polygon/convex-poly/angles/Sep06/\jobname.tex}

\begin{abstract}
Let $\P$ be a cyclic $n$-gon with $n\ge3$, the central angles  $\theta_0\in(-\pi,\pi],\dots,\theta_{n-1}\in(-\pi,\pi]$, and the winding number $w:=(\theta_0+\dots+\theta_{n-1})/(2\pi)$.
The vertices of $\P$ are assumed to be all distinct from one another. 
It is then proved that $\P$ is convex if and only if one of the following four conditions holds:
\begin{description}
\item[(I)] 
$w=1$ and $\th_0,\dots,\th_{n-1}>0$;
\item[(II)] 
$w=-1$ and $\th_0,\dots,\th_{n-1}<0$;
\item[(III)] 
$w=0$ and exactly one of the angles $\th_0,\dots,\th_{n-1}$ is negative;
\item[(IV)] 
$w=0$ and exactly one of the angles $\th_0,\dots,\th_{n-1}$ is positive.
\end{description}
\end{abstract}



\maketitle

\setcounter{section}{-1}

\section{Introduction}

The paper is organized as follows. 
In Section \ref{results}, the basic definitions are given and the main result, Theorem \ref{th:}, is stated.

In Section \ref{proofs}, the proofs are given. 
More specifically, Subsection \ref{subsec:lemmas} contains
statements of lemmas and, based on them, the proof of Theorem \ref{th:}; the proofs of all lemmas are deferred further to Subsection \ref{subsec:proofs of lemmas}.  

\section{Definitions and Results}\label{results}

A {\em polygon} is any finite sequence $\P:=(V_0,\dots,V_{n-1})$ of points (or, interchangeably, vectors) on the Euclidean plane. 
The points $V_0,\dots,V_{n-1}$ are called the {\em vertices} of $\P$.
The smallest value that we shall allow here for the integer $n$ is $3$. 
The segments, or closed intervals,
$$[V_i,V_{i+1}]:=\conv\{V_i,V_{i+1}\}\quad\text{for}\ i\in\{0,\dots,n-1\}$$
are called the {\em edges} of polygon $\P$, where 
$$V_n:=V_0.$$
The symbol $\conv$ denotes, as usual, the convex hull \cite[page 12]{rock}.
Note that, if $V_i=V_{i+1}$, then the edge $[V_i,V_{i+1}]$ is a singleton set. 
 
Our terminology concerning convexity corresponds to that in \cite{rock}. 
Here and in the sequel, we also use the notation
$$\intr km:=\{i\in\Z\colon k\le i\le m\},$$
where $\Z$ is the set of all integers; in particular, $\intr km$ is empty if $m<k$. 

Let us define the convex hull of polygon $\P$ as the convex hull of the set of its vertices: $\conv\P:=\conv\{V_0,\dots,V_{n-1}\}$.

Given the above notion of the polygon, a {\em convex polygon} can be defined as a polygon $\P$ such that the union of the edges of $\P$ coincides with the boundary $\partial\conv\P$ of the convex hull $\conv\P$ of $\P$; cf. e.g. \cite[page 5]{yaglom}. 
Thus, one has

\begin{definition}\label{def:conv}
A polygon $\P=(V_0,\dots,V_{n-1})$ is {\em convex} if 
$$\bigcup_{i\in\intr0{n-1}}[V_i,V_{i+1}]=\partial\conv\P.$$
\end{definition}

Let us emphasize that a polygon in this paper is a sequence and therefore ordered. In particular, even if all the vertices $V_0,\dots,V_{n-1}$ of a polygon $\P=(V_0,\dots,V_{n-1})$ are the extreme points of the convex hull of $\P$, it does not necessarily follow that $\P$ is convex. For example, consider the points $V_0=(0,0)$, $V_1=(1,0)$, $V_2=(1,1)$, and $V_3=(0,1)$ in $\R^2$. Then polygon $(V_0,V_1,V_2,V_3)$ is convex, while polygon $(V_0,V_2,V_1,V_3)$ is not.

\begin{definition}\label{def:ordinary}
Let us say that a polygon $\P=(V_0,\dots,V_{n-1})$ is {\em ordinary} if 
its vertices are all distinct from one another: $(i\ne j\ \&\ i\in\intr0{n-1}\ \&\ j\in\intr0{n-1}\,)\implies V_i\ne V_j$.
\end{definition}

\begin{remark}\label{rem:ordinary}
The set of edges of any polygon can be represented as the union of the sets of edges of ordinary polygons.
\end{remark}

\begin{definition}\label{def:cyclic}
A polygon $\P=(V_0,\dots,V_{n-1})$ is {\em cyclic} if  
all of its vertices, $V_0,\dots,V_{n-1}$, lie on a circle (of a stricly positive radius). 
\end{definition}

Let us now define the central angles of a cyclic polygon $\P=(V_0,\dots,V_{n-1})$. 
We shall identify the Euclidean plane with $\R^2$ in such a way that 
$$V_0=(1,0)$$
and the center of the circumscribed circle is $(0,0)$, so that the circumscribed circle of the cyclic polygon is the unit one.
Then \cite[Chapter VIII, \S2]{bour}, for each $i\in\intr0n$, there is some $\si_i\in\R$ such that
\begin{equation}\label{eq:V,si}
V_i=(\cos\si_i,\sin\si_i);
\end{equation}
in fact, $\si_i$ is determined up to an arbitrary additive integer multiple of $2\pi$. Thus, for any $i\in\intr0n$ and any semi-open interval $I$ of length $2\pi$, there is a uniquely determined number $\si_i$ in interval $I$ such that \eqref{eq:V,si} holds.

This allows one to define the $\si_i$'s in a unique way. Namely, set
\begin{equation*}\label{eq:si0}
\si_0:=0;
\end{equation*}
then, successively for every $i\in\intr0{n-1}$, given a value of $\si_i$, take $\si_{i+1}$ to be the uniquely determined number in the semi-open interval $(\si_i-\pi,\si_i+\pi]$ such that $V_{i+1}=(\cos\si_{i+1},\sin\si_{i+1})$. 

In particular, the winding number of polygon $\P$
$$w:=\frac{\si_n}{2\pi}$$
is then uniquely determined.
Note that $w$ must be an integer (because $V_n=V_0$, so that $\frac{\si_n}{2\pi}=\frac{\si_n-\si_0}{2\pi}\in\Z$). 

Now, for every $i\in\intr0{n-1}$, introduce (the radian measure of) the {\em central angle} $\th_i$ corresponding to the edge $[V_i,V_{i+1}]$ of the cyclic polygon $\P=(V_0,\dots,V_{n-1})$ by
\begin{equation*}\label{eq:th}
\th_i:=\si_{i+1}-\si_i;
\end{equation*}
then, by the construction of the $\si_i$'s,
\begin{equation}\label{eq:th-in}
\th_i\in(-\pi,\pi].
\end{equation}

Note also that the ordinariness of $\P$ will imply, in particular, that
$$\th_i\ne0\quad\forall i\in\intr0{n-1}.$$

Now we are prepared to state the main result of this paper, which was needed in \cite{jgeom}.

\begin{theorem}\label{th:}
Let $\P=(V_0,\dots,V_{n-1})$ be an ordinary cyclic polygon.
Then $\P$ is convex if and only if one of the following four conditions holds:
\begin{description}
\item[(I)] 
$w=1$ and $\th_0,\dots,\th_{n-1}>0$;
\item[(II)] 
$w=-1$ and $\th_0,\dots,\th_{n-1}<0$;
\item[(III)] 
$w=0$ and exactly one of the angles $\th_0,\dots,\th_{n-1}$ is negative:
$$\exists i\in\intr0{n-1}\quad 
(\th_i<0\ \&\ \forall j\in\intr0{n-1}\setminus\{i\}\ \th_j>0);$$
\item[(IV)] 
$w=0$ and exactly one of the angles $\th_0,\dots,\th_{n-1}$ is positive:
$$\exists i\in\intr0{n-1}\quad 
(\th_i>0\ \&\ \forall j\in\intr0{n-1}\setminus\{i\}\ \th_j<0).$$
\end{description}
\end{theorem}

Note that alternatives (III) and (IV) were overlooked in the known heuristic proof of \cite[Theorem~1 there]{macnab}. This 
suggests that there likely are genuine and substantial difficulties with the proof (and even formulation) of Theorem \ref{th:} that need to be overcome. 
It appears that at the root of these difficulties is the necessity to bridge the gap between such apparently distant notions as the polygon convexity and the central angles.  
Moreover, papers \cite{elim} and \cite{test} suggest that the notion of polygon convexity is rather complex by itself, as it connects the notion of a polygon (and hence that of order) with the notion of convexity. 

\section{Proofs}\label{proofs}

\subsection{Lemmas, and the Proof of Theorem \ref{th:} } \label{subsec:lemmas}

Introduce the determinants
\begin{equation}\label{eq:De}
\De_{\al,\be,\g}:=
\left|
\begin{matrix}
1&\cos\si_\al&\sin\si_\al \\
1&\cos\si_\be&\sin\si_\be \\
1&\cos\si_\g&\sin\si_\g 
\end{matrix}\right|
\end{equation}
for $\al$, $\be$, and $\g$ in the set $\intr0{n-1}$. 

\begin{lemma}\label{lem:De}
For $i$ and $j$ in $\intr0{n-1}$,
$$\De_{j,i,i\oplus1}=
4 \sin\frac{\th_i}2\, \sin\frac{\si_j-\si_i}2\, \sin\frac{\si_j-\si_{i+1}}2,
$$
where
$$i\oplus1:=
\begin{cases}
i+1 &\text{if}\ i\in\intr0{n-2},\\
0 &\text{if}\ i=n-1. 
\end{cases}
$$ 
\end{lemma}

\begin{lemma}\label{lem:if}
The determinants 
$$\De_{j,i,i\oplus1}\quad\text{for}\quad 
i\in\intr0{n-1}\quad\text{and}\quad j\in\intr0{n-1}\setminus\{i,i\oplus1\}$$
are all strictly positive under condition (I) or (III), and these determinants are all strictly negative under condition (II) or (IV). 
\end{lemma}

\begin{lemma}\label{lem:<0}
Suppose that, for some $i$, $j$, and $k$ in the set $\intr0{n-1}$, there exist integers $p$ and $q$ such that
\begin{align} 
\hat\si_j:=\si_j+2\pi p\in(\si_i,\si_{i+1})\quad & \&\quad 
\hat\si_k:=\si_k+2\pi q\in(\si_{i+1},\si_i+2\pi) \label{eq:or1}\\
\intertext{or}
\hat\si_j\in(\si_{i+1},\si_i)\quad & \&\quad 
\hat\si_k\in(\si_i,\si_{i+1}+2\pi). \label{eq:or2}
\end{align}
Then
$$\De_{j,i,i\oplus1}\,\De_{k,i,i\oplus1}<0.$$ 
\end{lemma}

\begin{lemma}\label{lem:patterns}
Suppose that one has one of the following three patterns for the $\th_i$'s:
\begin{description}
\item[(P1)]
$\th_0,\dots,\th_m>0$ and $\th_0+\dots+\th_m>2\pi$, for some $m\in\intr0{n-1}$;
\item[(P2)]
$\th_0<0$; $\th_1,\dots,\th_m>0$; and $\th_{m+1}<0$, for some $m\in\intr2{n-2}$;
\item[(P3)]
$\th_0<0$, $\th_1>0$, $\th_2<0$, and $\th_3>0$ (so that $n\ge4$).
\end{description}
Then either \eqref{eq:or1} or \eqref{eq:or2} holds for some $i$, $j$, and $k$ in the set $\intr0{n-1}$ and some integers $p$ and $q$.
\end{lemma}

\begin{definition}\label{def:to one side}
Let $P_0,\dots,P_m$ be distinct points on the plane. Let us write $P_2,\dots,P_m\s[P_0,P_1]$ 
and say
that
points $P_2,\dots,P_m$ are {\em strictly to one side}) 
of segment $[P_0,P_1]$ if 
the (straight) line through points $P_0$ and $P_1$ is the boundary of an open 
half-plane containing the set $\{P_2,\dots,P_m\}$. 

For any given $i\in\intr0{n-1}$, let us say that
a polygon $\P=(V_0,\dots,V_{n-1})$ is strictly to one side 
of its edge $[V_i,V_{i+1}]$ if the set 
$\{V_j\colon j\in\intr0{n-1}\setminus\{i,i\oplus1\}\}$ 
is so.

Let us say that a polygon is {\em strictly to-one-side} if it is strictly to one side of every one of its edges. 
\end{definition}

\begin{lemma}\label{lem:to one side}
An ordinary cyclic polygon $\P=(V_0,\dots,V_{n-1})$ 
is convex if and only if it is strictly to one side.    
\end{lemma}

\begin{lemma}\label{prop:calculation}
Let $x_i$ and $y_i$ denote the coordinates of points $V_i$, so that
$V_i=(x_i,y_i)$ for all $i\in\intr0{n-1}$.
Then, for any choice of $\al$, $\be$, $i$, and $j$ in $\intr0{n-1}$, 
$$V_\al,V_\be\s[V_i,V_j]\iff
\De_{\al,i,j}\,\De_{\be,i,j}>0,$$
where $\De_{\al,i,j}$ are given by \eqref{eq:De}.
\end{lemma}

\begin{lemma}\label{lem:calculation}
An ordinary cyclic polygon $\P=(V_0,\dots,V_{n-1})$ 
is convex if and only if for each $i\in\intr0{n-1}$ the signs of the determinants $\De_{j,i,i\oplus1}$ are the same for all $j\in\intr0{n-1}\setminus\{i,i\oplus1\}$.   
\end{lemma}

\begin{lemma}\label{lem:necess}
An ordinary cyclic polygon $\P=(V_0,\dots,V_{n-1})$ is convex {\em only} if none of the three patterns, (P1)--(P3), listed in Lemma \ref{lem:patterns} takes place.
\end{lemma}

\begin{proof}[Proof of Theorem \ref{th:}]

\textbf{``If"}\quad The ``if" part of Theorem \ref{th:} follows immediately from Lemmas \ref{lem:if} and \ref{lem:calculation}.

\textbf{``Only if"}\quad To prove the ``only if" part of Theorem \ref{th:}
suppose, to the contrary, that $\P=(V_0,\dots,V_{n-1})$ is an ordinary cyclic polygon with $n\ge3$, but none of the conditions (I)--(IV) holds. Then at least one of the following 5 alternatives must take place:
\begin{description}
\item[(A1)]
$\th_0,\dots,\th_{n-1}>0$ and $w\ne1$;
\item[(A2)]
$\th_0,\dots,\th_{n-1}<0$ and $w\ne-1$;
\item[(A3)]
among the $n$ values $\th_0,\dots,\th_{n-1}$ exactly one value is negative, and $w\ne0$;
\item[(A4)]
among the $n$ values $\th_0,\dots,\th_{n-1}$ exactly one value is positive, and $w\ne0$;
\item[(A5)]
among the $n$ values $\th_0,\dots,\th_{n-1}$ at least two values are positive and at least two values are negative (so that $n\ge4$).
\end{description}
It suffices to show that each of these 5 alternatives leads to a contradiction.

\textbf{Alternative (A1):}\quad Assume that (A1) takes place. Then $w\in\{2,3,\dots\}$, so that 
$\th_0+\dots+\th_{n-1}>2\pi$. Thus, pattern (P1) listed in Lemma \ref{lem:patterns} takes place (with $m=n-1$). 
By Lemma \ref{lem:necess}, this contradicts the convexity of polygon $\P$. 

\textbf{Alternative (A2):}\quad This is quite similar to alternative (A1). In fact, (A2) reduces to (A1) if polygon $\P=(V_0,\dots,V_{n-1})$ is replaced by the ``re-oriented" polygon $(V_{n-1},\dots,V_0)$. 

\textbf{Alternative (A3):}\quad Assume that (A3) takes place. Then, w.l.o.g., $\th_0,\dots,\th_{n-2}>0$;  $\th_{n-1}<0$; $\si_n=\th_0+\dots+\th_{n-1}=2\pi w$, where $w$ is a {\em nonzero} integer. Moreover, 
$2\pi w=\si_n=\si_{n-1}+\th_{n-1}>\th_{n-1}>-\pi$, so that $w$ is a nonnegative integer. Then in fact $w\ge1$ (since 
$w\ne0$). 
Thus, 
$$\th_0+\dots+\th_{n-2}=\si_{n-1}=\si_n-\th_{n-1}>\si_n=2\pi w\ge2\pi,$$
so that pattern (P1) listed in Lemma \ref{lem:patterns} takes place (with $m=n-2$). By Lemma \ref{lem:necess}, this contradicts the convexity of polygon $\P$. 

\textbf{Alternative (A4):}\quad This is quite similar to alternative (A3). In fact, (A4) reduces to (A3) if polygon $\P=(V_0,\dots,V_{n-1})$ is replaced by the ``re-oriented" polygon $(V_{n-1},\dots,V_0)$. 

\textbf{Alternative (A5):}\quad Assume that (A5) takes place. For the given convex cyclic polygon $\P=(V_0,\dots,V_{n-1})$, let us call a non-empty set $R\subseteq\intr0{n-1}$ a {\em plus-run} if $\th_i>0$ for all $i\in R$ and $R=\intr k\ell$ for some integers $k$ and $\ell$. Let us refer to the cardinality of a plus-run as its {\em length}. A plus-run will be called {\em maximal} if it is not contained in any other plus-run. Similarly defined are a {\em minus-run} (with condition $\th_i<0$ in place of $\th_i>0$), its length, and a maximal minus-run. A {\em run} is a set which is either a plus-run or a minus-run. A {\em maximal run} is a set which is either a maximal plus-run or a maximal minus-run. Let $r$ denote the number of maximal runs for the polygon $\P$. 

Note that (A5) implies that $r\ge2$. Therefore, it suffices to consider the following three cases.

\textbf{Case (A5-1): There are no runs of length $\ge2$.}\quad That is, the signs of the $\th_i$'s alternate. In view of the ``re-orientation" possibility, one may assume w.l.o.g. that pattern (P3) listed in Lemma \ref{lem:patterns} takes place. By Lemma \ref{lem:necess}, this contradicts the convexity of polygon $\P$. 

\textbf{Case (A5-2): $r=2$.}\quad In view of the ``re-orientation" possibility, one may assume w.l.o.g. that $\th_0,\dots,\th_{k-1}<0$ and $\th_k,\dots,\th_{n-1}>0$, for some $k\in\intr2{n-2}$. 
Consider now the cyclic permutation of the $\th_i$'s:
$$(\hat\th_0,\hat\th_1,\dots,\hat\th_{n-k},\hat\th_{n-k+1},\dots,\hat\th_{n-1})
:=
(\th_{k-1},\th_k,\dots,\th_{n-1},\th_0,\dots,\th_{k-2}),$$
corresponding to the cyclic permutation 
$$\hat\P:=(\hat V_0,\hat V_1,\dots,\hat V_{n-k},\hat V_{n-k+1},\dots,\hat V_{n-1})
:=
(V_{k-1},V_k,\dots,V_{n-1},V_0,\dots,V_{k-2})$$
of polygon $\P=(V_0,\dots,V_{n-1})$.
Then $\hat\th_0<0$; $\hat\th_1,\dots,\hat\th_{n-k}>0$; and $\hat\th_{n-k+1}<0$.
That is, pattern (P2) listed in Lemma \ref{lem:patterns} takes place for the cyclic polygon $\hat\P$ with $m:=n-k$ (note that $m\in\intr2{n-2}$, since $k\in\intr2{n-2}$). By Lemma \ref{lem:necess}, this is a contradiction, since the convexity of polygon $\P$ implies that of polygon $\hat\P$.
 
\textbf{Case (A5-3): $r\ge3$ and there is a (maximal) run of length $\ge2$.}\quad 
Let 
$$\{0,\dots,i_1\},\ \{i_1+1,\dots,i_2\},\ \{i_2+1,\dots,i_3\},\ \dots,\ \{i_{r-1}+1,\dots,n-1\}$$
be the maximal runs, for some $i_1,i_2,i_3,\dots,i_{r-1}$ in $\intr0{n-1}$ such that $0\le i_1<i_2<i_3<\dots<i_{r-1}<n-1$. In view of the possibility of a cyclic permutation, one may assume w.l.o.g. that the second maximal run is of length $\ge2$; that is, $i_2\ge i_1+2$. In view of the ``re-orientation" possibility, one may assume w.l.o.g. that $\th_{i_1}<0$; $\th_{i_1+1},\dots,\th_{i_2}>0$; and $\th_{i_2+1}<0$. 
Thus, pattern (P2) listed in Lemma \ref{lem:patterns} takes place for a cyclic permutation of polygon $\P$ (with $m=i_2-i_1\in\intr2{n-2}$). By Lemma \ref{lem:necess}, this contradicts the convexity of polygon $\P$. 
\end{proof}

\subsection{Proofs of the Lemmas}\label{subsec:proofs of lemmas}

\begin{proof}[Proof of Lemma \ref{lem:De}]
Subtracting the second row of the determinant $\De_{j,i,i\oplus1}$ from the other two ones and then expanding the determinant along the first column and using the product expressions for $\cos x-\cos y$ and $\sin x-\sin y$, one has
\begin{align*}
& \De_{j,i,i\oplus1}\\
=& 4\sin\frac{\th_i}2\,\sin\frac{\si_j-\si_i}2\, 
\lp\sin\frac{\si_i+\si_j}2\,\cos\frac{\si_i+\si_{i+1}}2\,
-\sin\frac{\si_i+\si_{i+1}}2\,\cos\frac{\si_i+\si_j}2\,\rp \\
=& 4 \sin\frac{\th_i}2\, \sin\frac{\si_j-\si_i}2\, \sin\frac{\si_j-\si_{i+1}}2,
\end{align*}
where we also used the identity  $\sin x\cos y-\sin y\cos x=\sin(x-y)$. 
\end{proof}

\begin{proof}[Proof of Lemma \ref{lem:if}] Take any $i$ and $j$ such that
$$i\in\intr0{n-1}\quad\text{and}\quad j\in\intr0{n-1}\setminus\{i,i\oplus1\}.$$

\textbf{(I):}\quad Consider the case when conditions (I) hold. Recall that condition $w=1$ means that $\th_0+\dots+\th_{n-1}=2\pi$. Consider next the two possible subcases: $i=n-1$ and $i\in\intr0{n-2}$. 

If $i=n-1$, then $i\oplus1=0$ and $j\in\intr1{n-2}=\intr1{i-1}$. Hence, in view of \eqref{eq:th-in} and Lemma \ref{lem:De}, 
\begin{align*}
\frac{\th_i}2&\in(0,\pi/2]\subset(0,\pi),\\
\frac{\si_j-\si_i}2&=-\frac12(\th_{j+1}+\dots+\th_i)\in(-\pi,0),\\
\frac{\si_j-\si_{i+1}}2&=-\frac12(\th_{j+1}+\dots+\th_{i+1})\in(-\pi,0),\\
\De_{j,i,i\oplus1}&>0.
\end{align*}

If $i\in\intr0{n-2}$, then $i\oplus1=i+1\in\intr1{n-1}$ and $j\in\intr0{i-1}\cup\intr{i+2}{n-1}$. 
If $j\in\intr0{i-1}$, then both $(\si_j-\si_i)/2$ and $(\si_j-\si_{i+1})/2$ lie in the interval $(-\pi,0)$.
If $j\in\intr{i+2}{n-1}$, then both $(\si_j-\si_i)/2$ and $(\si_j-\si_{i+1})/2$ lie in the interval $(0,\pi)$. 

It follows that $\De_{j,i,i\oplus1}>0$ if conditions (I) hold. 

\textbf{(II):}\quad The case when conditions (II) hold is quite similar to (I). In this case, one has $\De_{j,i,i\oplus1}<0$ for all $i\in\intr0{n-1}$ and $j\in\intr0{n-1}\setminus\{i,i\oplus1\}$. 

\textbf{(III):}\quad Consider the case when conditions (III) hold. 
Here, w.l.o.g., it is $\th_{n-1}$ that is negative, of all the central angles $\th_0,\dots,\th_{n-1}$. That is, one has $\th_0,\dots,\th_{n-2}>0$ and $\th_{n-1}<0$. Also, condition $w=0$ means that $\th_0+\dots+\th_{n-1}=0$. 
Hence, $\si_{n-1}=-\th_{n-1}\in(0,\pi)$, and so, $\si_s-\si_r\in(0,\pi)$ whenever $0\le r<s\le n-1$.  
Consider next the two possible subcases: $i=n-1$ and $i\in\intr0{n-2}$. 

If $i=n-1$, then $i\oplus1=0$ and $j\in\intr1{n-2}=\intr1{i-1}$. Hence, 
$\th_i=\th_{n-1}\in(-\pi,0)$, $\si_j-\si_i\in(-\pi,0)$, and $\si_j-\si_{i+1}=\si_j-\si_0\in(0,\pi)$. Thus, $\De_{j,i,i\oplus1}>0$. 

If $i\in\intr0{n-2}$, then $i\oplus1=i+1\in\intr1{n-1}$ and $j\in\intr0{i-1}\cup\intr{i+2}{n-1}$. 
If $j\in\intr0{i-1}$, then both $\si_j-\si_i$ and $\si_j-\si_{i+1}$ lie in the interval $(-\pi,0)$.
If $j\in\intr{i+2}{n-1}$, then both $\si_j-\si_i$ and $\si_j-\si_{i+1}$ lie in the interval $(0,\pi)$. 

It follows that $\De_{j,i,i\oplus1}>0$ if conditions (III) hold. 

\textbf{(IV):}\quad The case when conditions (IV) hold is quite similar to (III). In this case, one has $\De_{j,i,i\oplus1}<0$ for all $i\in\intr0{n-1}$ and $j\in\intr0{n-1}\setminus\{i,i\oplus1\}$.
\end{proof}

\begin{proof}[Proof of Lemma \ref{lem:<0}] 
Assume that the conditions of Lemma \ref{lem:<0} hold. 
Note that, in view of Lemma \ref{lem:De} and identity $\sin\frac{x+2\pi p}2=(-1)^p\sin\frac x2$ for $p\in\Z$, one has
$$\frac{\De_{j,i,i\oplus1}\,\De_{k,i,i\oplus1}}{16\sin^2\frac{\th_i}2}
=\sin\frac{\hat\si_j-\si_i}2\, \sin\frac{\hat\si_j-\si_{i+1}}2\, 
\sin\frac{\hat\si_k-\si_i}2\, \sin\frac{\hat\si_k-\si_{i+1}}2.$$ 

Consider now the following two possible cases.

\textbf{Case 1: \eqref{eq:or1} holds.} \quad Here condition $\hat\si_j\in(\si_i,\si_{i+1})$ implies that $\si_i<\si_{i+1}$. Moreover, one has
\begin{align*}
0<\hat\si_j-\si_i<\si_{i+1}-\si_i=\th_i\le\pi,\quad
&\text{so that}\quad 
\frac{\hat\si_j-\si_i}2\in(0,\pi/2); \\
0>\hat\si_j-\si_{i+1}>\si_i-\si_{i+1}=-\th_i\ge-\pi,\quad
&\text{so that}\quad 
\frac{\hat\si_j-\si_{i+1}}2\in(-\pi/2,0); \\
0<\si_{i+1}-\si_i<\hat\si_k-\si_i<2\pi,\quad
&\text{so that}\quad 
\frac{\hat\si_k-\si_i}2\in(0,\pi); \\
0<\hat\si_k-\si_{i+1}<\si_i-\si_{i+1}+2\pi<2\pi,\quad
&\text{so that}\quad 
\frac{\hat\si_k-\si_{i+1}}2\in(0,\pi).
\end{align*}
Thus, Lemma \ref{lem:<0} follows in Case 1.  

\textbf{Case 2: \eqref{eq:or2} holds.} \quad This case is quite similar to Case 1; interchange $\si_i$ with $\si_{i+1}$ and $\th_i$ with $(-\th_i)$ everywhere.  \end{proof}

\begin{proof}[Proof of Lemma \ref{lem:patterns}] 
\textbf{(P1):}\quad Suppose that pattern (P1) takes place:
$$\text{$\th_0,\dots,\th_m>0$ and $\th_0+\dots+\th_m>2\pi$, for some $m\in\intr0{n-1}$.}$$
Then, in view of \eqref{eq:th-in}, the condition $[\si_{m+1}=]\th_0+\dots+\th_m>2\pi$ implies that $m\ge2$. Moreover, w.l.o.g., $m$ is the smallest integer $\ge2$ such that (P1) holds. Hence,
$$\si_m\le2\pi,$$
and so, 
$\si_{m+1}=\si_m+\th_m\le2\pi+\pi=3\pi$; that is, 
$$\si_{m+1}-2\pi\in(0,\pi].$$
Also, $\si_0=0$, $\si_1=\th_0\in(0,\pi]$, and $\si_m=\si_{m+1}-\th_m>2\pi-\th_m\ge\pi$.
Hence, 
$$(0,\pi]\subseteq(0,\si_m]=(\si_0,\si_m]
=\bigcup_{i=0}^{m-1}(\si_i,\si_{i+1}].$$
Since $\si_{m+1}-2\pi\in(0,\pi]$ and $\si_{m+1}-2\pi\ne\si_i$ for any $i\in\intr0{m-1}$ (because polygon $\P$ is ordinary), 
it follows that
$$\si_{m+1}-2\pi\in(\si_i,\si_{i+1})$$
for some 
$$i\in\intr0{m-1};$$
that is, the first half of condition \eqref{eq:or1} takes place with $j=m+1$ and $p=-1$. 

Let us show that the second half of condition \eqref{eq:or1} takes place for appropriate $k$ and $q$. Here we must distinguish between the two possible cases: $i\in\intr0{m-2}$ and $i=m-1$.

\textbf{Case (P1-1):} $i\in\intr0{m-2}$.\quad Here, 
$$\si_{i+1}\le\si_m\le2\pi\le\si_i+2\pi.$$
Therefore, $\si_m\in[\si_{i+1},\si_i+2\pi]$, and so, $\si_m\in(\si_{i+1},\si_i+2\pi)$, because $i+1<m$ and polygon $\P$ is ordinary. Thus, one has the second half of condition \eqref{eq:or1} with $k=m$ and $q=0$.     

\textbf{Case (P1-2):} $i=m-1$.\quad In view of \eqref{eq:th-in}, one has $\th_{m-2},\th_{m-1}\le\pi$; moreover, conditions $\th_{m-2}=\th_{m-1}=\pi$ would imply that $V_{m-2}=V_m$, which would contradict the ordinariness of polygon $\P$. Therefore,  
$\th_{m-2}+\th_{m-1}<2\pi$. Hence,
$$\si_{m-2}+2\pi>\si_{m-2}+\th_{m-2}+\th_{m-1}=\si_m,$$
and so, $\si_{m-2}+2\pi\in(\si_m,\si_{m-1}+2\pi)$, which implies the second half of condition \eqref{eq:or1} with $k=m-2$ and $q=1$. 

This completes the consideration of pattern (P1).

\textbf{(P2):}\quad Suppose that pattern (P2) takes place:
$$\text{$\th_0<0$; $\th_1,\dots,\th_m>0$; and $\th_{m+1}<0$, for some $m\in\intr2{n-2}$.}$$
Here we must distinguish the following four possible cases.

\textbf{Case (P2-1):} $\si_2>0$.\quad Here, $\si_0=0\in(\si_1,\si_2)$ (since $\si_1=\th_0<0$ and $\si_2>0$); that is, the first half of condition \eqref{eq:or1} takes place with $i=1$, $j=0$, and $p=0$. 
Next, note that $\th_1+\th_2<2\pi$ (because $\th_1,\th_2\le\pi$ and $\th_1=\th_2=\pi$ would imply $V_2=V_0$, which would contradict the ordinariness of $\P$). Hence,
$\si_3=\si_1+\th_1+\th_2<\si_1+2\pi$.
Also, $\si_3-\si_2=\th_2>0$. 
It follows that $\si_3\in(\si_2,\si_1+2\pi)$; that is, the second half of condition \eqref{eq:or1} takes place with $k=3$ and $q=0$ (and the same $i=1$). 

\textbf{Case (P2-2):} $\si_{m+1}<0$, $\si_{m+2}>\si_1$.\quad Then 
$$\bigcup_{i=1}^m(\si_i,\si_{i+1}]=(\si_1,\si_{m+1}]\supseteq(\si_1,\si_{m+1})\ni\si_{m+2},$$
because $\si_{m+2}=\si_{m+1}+\th_{m+1}<\si_{m+1}$ and by condition $\si_{m+2}>\si_1$.
Hence and because of the ordinariness of $\P$, one has $$\si_{m+2}\in(\si_i,\si_{i+1})$$
for some $i\in\intr1m$. 
On the other hand, $\si_{i+1}=\si_{m+1}-(\th_{i+1}+\dots+\th_m)\le\si_{m+1}<0$, by the definition of Case (P2-2). Also,
$$\si_i+2\pi=\si_1+\th_1+\dots+\th_{i-1}+2\pi\ge\si_1+2\pi
=\th_0+2\pi>-\pi+2\pi>0,$$
so that 
$$\si_0=0\in(\si_{i+1},\si_i+2\pi).$$ 
Thus, condition \eqref{eq:or1} takes place with $j=m+2$, $p=0$, $k=0$, and $q=0$. 

\textbf{Case (P2-3):} $\si_{m+1}<0$, $\si_{m+2}<\si_1$.\quad Then 
$\si_{m+1}=\si_1+\th_1+\dots+\th_m>\si_1$ and $\si_{m+1}<0=\si_0$, so that 
$$\si_{m+1}\in(\si_1,\si_0).$$ 
On the other hand, $\si_{m+2}=\si_{m+1}+\th_{m+1}>\si_1+\th_{m+1}>\si_1-\pi=\th_0-\pi>-2\pi
=\si_0-2\pi$ and $\si_{m+2}<\si_1$ (by the definition of Case (P2-3)); hence,
$$\si_{m+2}+2\pi\in(\si_0,\si_1+2\pi).$$ 
Thus, condition \eqref{eq:or2} takes place with $i=0$, $j=m+1$, $p=0$, $k=m+2$, and $q=1$. 

\textbf{Case (P2-4):} $\si_2<0$, $\si_{m+1}>0$.\quad Let here $i$ be the greatest integer such that $\si_i<0$ (such an $i$ exists, by the definition of Case (P2-4)). Then $\si_{i+1}>0$ and $i\in\intr2m$. Hence, $$\si_0=0\in(\si_i,\si_{i+1}).$$ 
On the other hand, $\si_1=\th_0<\th_0+\th_1+\dots+\th_{i-1}=\si_i$ and 
$\si_1=\th_0\ge-\pi=0+\pi-2\pi>\si_i+\th_i-2\pi=\si_{i+1}-2\pi$, so that 
$$\si_1+2\pi\in(\si_{i+1},\si_i+2\pi).$$
Thus, condition \eqref{eq:or1} takes place with $j=0$, $p=0$, $k=1$, and $q=1$. 

This completes the consideration of pattern (P2).

\textbf{(P3):}\quad Suppose that pattern (P3) takes place:
$$\text{$\th_0<0$, $\th_1>0$, $\th_2<0$, and $\th_3>0$ (so that $n\ge4$).}$$
Here we must distinguish the following 6 possible cases.

\textbf{Case (P3-1):} $\si_3>0$, $\si_4>\si_2$.\quad Here, 
$\si_2=\si_3-\th_2>\si_3>0=\si_0$. Also, $\si_1=\th_0<0=\si_0$. Hence,
$$\si_0\in(\si_1,\si_2).$$
On the other hand, $\si_4>\si_2$ (by the definition of Case (P3-1)) and 
$\si_4=\si_1+\th_1+\th_2+\th_3<\si_1+\th_1+\th_3\le\si_1+2\pi$. It follows that
$$\si_4\in(\si_2,\si_1+2\pi).$$  
Thus, condition \eqref{eq:or1} takes place with $i=1$, $j=0$, $p=0$, $k=4$, and $q=0$.

\textbf{Case (P3-2):} $\si_3>0$, $\si_4<\si_2$.\quad Here, 
$\si_4=\si_3+\th_3>\si_3$ and, by the definition of Case (P3-2), $\si_4<\si_2$. Hence,
$$\si_4\in(\si_3,\si_2).$$
On the other hand, $\si_2-2\pi=\th_0+\th_1-2\pi<\th_1-2\pi\le\pi-2\pi<0=\si_0$ and, by the definition of Case (P3-2), $\si_3>0=\si_0$. 
It follows that
$$\si_0+2\pi\in(\si_2,\si_3+2\pi).$$  
Thus, condition \eqref{eq:or2} takes place with $i=2$, $j=4$, $p=0$, $k=0$, and $q=1$.

\textbf{Case (P3-3):} $\si_2<0$, $\si_3<\si_1$.\quad Here, 
$\si_2=\si_1+\th_1>\si_1$ and, by the definition of Case (P3-3), $\si_2<0=\si_0$. Hence,
$$\si_2\in(\si_1,\si_0).$$
On the other hand, $\si_3=\si_0+\th_0+\th_1+\th_2>\si_0+\th_0+\th_2>\si_0-2\pi$ and, by the definition of Case (P3-3), $\si_3<\si_1$. 
It follows that
$$\si_3+2\pi\in(\si_0,\si_1+2\pi).$$  
Thus, condition \eqref{eq:or2} takes place with $i=0$, $j=2$, $p=0$, $k=3$, and $q=1$.

\textbf{Case (P3-4):} $\si_2<0$, $\si_3>\si_1$.\quad Here, 
$\si_3=\si_2+\th_2<\si_2$ and, by the definition of Case (P3-4), $\si_3>\si_1$. Hence,
$$\si_3\in(\si_1,\si_2).$$
On the other hand, $\si_1+2\pi=\si_0+\th_0+2\pi>\si_0-\pi+2\pi>\si_0$ and, by the definition of Case (P3-4), $\si_2<0=\si_0$. 
It follows that
$$\si_0\in(\si_2,\si_1+2\pi).$$  
Thus, condition \eqref{eq:or1} takes place with $i=1$, $j=3$, $p=0$, $k=0$, and $q=0$.

\textbf{Case (P3-5):} $\si_2>0$, $\si_3<0$, $\si_3>\si_1$.\quad Here, by the definition of Case (P3-5), $\si_3<0=\si_0$ and $\si_3>\si_1$. Hence,
$$\si_3\in(\si_1,\si_0).$$
On the other hand, $\si_2=\si_1+\th_1\le\si_1+\pi<\si_1+2\pi$ and, by the definition of Case (P3-5), $\si_2>0=\si_0$. 
It follows that
$$\si_2\in(\si_0,\si_1+2\pi).$$  
Thus, condition \eqref{eq:or2} takes place with $i=0$, $j=3$, $p=0$, $k=2$, and $q=0$.

\textbf{Case (P3-6):} $\si_2>0$, $\si_3<0$, $\si_3<\si_1$.\quad Here, 
$\si_0=0>\th_0=\si_1$ and,
by the definition of Case (P3-6), $\si_2>0=\si_0$. Hence,
$$\si_0\in(\si_1,\si_2).$$
On the other hand, $\si_3=\si_2+\th_2>\si_2-\pi>\si_2-2\pi$ and, by the definition of Case (P3-6), $\si_3<\si_1$. 
It follows that
$$\si_3+2\pi\in(\si_2,\si_1+2\pi).$$  
Thus, condition \eqref{eq:or1} takes place with $i=1$, $j=0$, $p=0$, $k=3$, and $q=1$.

This completes the consideration of pattern (P3) as well.
\end{proof}

\begin{proof}[Proof of Lemma \ref{lem:to one side}]

{\bf ``Only if":}\quad Suppose a polygon $\P=(V_0,\dots,V_{n-1})$ is ordinary and convex. Take any $i\in\intr0{n-1}$. Then, by Definition \ref{def:conv}, $[V_i,V_{i+1}]\subseteq\partial\conv\P$. Hence, by \cite[Theorem 11.6]{rock}, the line $\ell$ through points $V_i$ and $V_{i+1}$ is the boundary of a closed half-plane containing
$\conv\P$. Moreover, for any $j\in\intr0{n-1}\setminus\{i,i\oplus1\}$, point $V_j$ is not on line $\ell$ (because (i) no line can have more than two distinct points in common with a circle and
(ii) polygon $\P$ is cyclic and ordinary). 
Hence, polygon $\P$ is strictly to one side of its edge $[V_i,V_{i+1}]$, for each $i\in\intr0{n-1}$.

{\bf ``If":}\quad Suppose a polygon $\P=(V_0,\dots,V_{n-1})$ is strictly to-one-side. For any $i\in\intr0{n-1}$, consider the line $\ell$ through points $V_i$ and $V_{i+1}$. Then $\ell$ is the boundary of a closed half-plane $H$ containing
$\conv\P$. By the definition of the convex hull, $[V_i,V_{i+1}]\subseteq\conv\P$. Hence,
$$[V_i,V_{i+1}]\subseteq\ell\cap\conv\P=\partial H\cap\conv\P\subseteq\partial \conv\P$$
(the latter inclusion follows because $\conv\P\subseteq H$). 
\end{proof}

\begin{proof}[Proof of Lemma \ref{prop:calculation}]
Take any $\al$, $\be$, $i$, $j$ in the set $\intr0{n-1}$. By Definition \ref{def:to one side}, one has $V_\al,V_\be\s[V_i,V_j]$ if and only if $V_j\ne V_i$ and there exists some vector $\nn=(a,b)\in\R^2$ such that
$$\nn\cdot\v ij=0<\nn\cdot\v i\g\quad\text{for}\ \g\in\{\al,\be\}.$$
Since
$$\De_{\al,i,j}=
\left|
\begin{matrix}
1&x_\al-x_i&y_\al-y_i \\
1&0&0 \\
1&x_j-x_i&y_j-y_i 
\end{matrix}\right|,
$$
one may replace w.l.o.g. the points $V_\al$, $V_\be$, $V_i$, $V_j$ by $V_\al-V_i$, $V_\be-V_i$, $V_i-V_i=(0,0)$, $V_j-V_i$, respectively. Hence, w.l.o.g., 
$$V_i=(0,0).$$
Then the condition $V_\al,V_\be\s[V_i,V_j]$ can be rewritten as 
$$ax_j+by_j=0<ax_\g+by_\g\quad\text{for}\ \g\in\{\al,\be\}$$
and $(x_j,y_j)\ne(0,0)$. W.l.o.g., $y_j\ne0$. 
Then condition $ax_j+by_j=0$ is equivalent to $b=-\frac{x_j}{y_j}a$, so that the inequality $0<ax_\g+by_\g$ can be rewritten as 
$\frac a{y_j}(x_\g y_j-x_j y_\g)>0$, or as $\frac a{y_j}\De_{\g,i,j}<0$ (where $\g\in\{\al,\be\}$); in particular, it follows that $a\ne0$.

We see that the condition $V_\al,V_\be\s[V_i,V_j]$ implies 
$$\De_{\al,i,j}\De_{\be,i,j}
=\lp\frac{y_j}a\rp^2 \lp\frac a{y_j}\De_{\al,i,j}\rp\lp\frac a{y_j}\De_{\be,i,j}\rp
>0.$$
This proves the ``$\Longrightarrow$" part of Lemma \ref{prop:calculation}. 

To prove the ``$\Longleftarrow$" part, let $\nn:=\vp(-y_j,x_j)$, where $\vp:=\sign\De_{\al,i,j}$. Then the condition $\De_{\al,i,j}\De_{\be,i,j}>0$ implies that $\vp=\sign\De_{\be,i,j}$. Also, $\nn\cdot\v ij=0$, while
$$\nn\cdot\v i\g=\vp(x_j y_\g-y_j x_\g)=\vp\De_{\g,i,j}=|\De_{\g,i,j}|>0$$
for $\g\in\{\al,\be\}$, so that the condition $V_\al,V_\be\s[V_i,V_j]$ takes place.
\end{proof}

\begin{proof}[Proof of Lemma \ref{lem:calculation}]
This follows immediately from Lemmas \ref{lem:to one side} and \ref{prop:calculation}. 
\end{proof}
 
\begin{proof}[Proof of Lemma \ref{lem:necess}]
This follows immediately from Lemmas \ref{lem:patterns}, \ref{lem:<0}, and \ref{lem:calculation}. 
\end{proof}

\renewcommand{\refname}{\textsf{\bf Literature}}

\bigskip

{\parskip0pt \parindent0pt \it Department of Mathematical Sciences

Michigan Technological University

Houghton, MI 49931

USA

e-mail: ipinelis@mtu.edu}

\end{document}